\documentclass[final,leqno]{siamltex}
\usepackage{amsfonts,amssymb,amsmath,cmmib57,color,graphics,mathdots}
\usepackage{amstext}
\usepackage{amscd}
\usepackage{graphicx}
\usepackage{longtable}
\usepackage{calc}
\usepackage[latin1]{inputenc}
\usepackage{url}
\usepackage{diagbox}
\usepackage{hyperref}
\usepackage{floatrow}
\usepackage{layout}
\usepackage{enumitem}

\usepackage{tikz}
\usetikzlibrary{decorations.pathreplacing}
\newcommand{\tikznode}[3][inner sep=0pt]{\tikz[remember
picture,baseline=(#2.base)]{\node(#2)[#1]{$#3$};}}

\usepackage{ytableau}

\usepackage{tikz}
\usetikzlibrary{tikzmark}

\newcommand{\hide}[1]{}

\newtheorem{example}{Example}

\newsavebox{\savepar}

\newcommand{\CC}{\mathbb{C}}

\newcommand{\cC}{{\cal C}}

\newcommand{\cB}{{\cal B}}

\newcommand{\la}{\lambda}

\newcommand{\orb}{{\cal O}}
\newcommand{\codim}{\textnormal{codim }}

\newcounter{algo}[section]

\usepackage{todonotes}

\title{On bundle closures of matrix pencils and matrix polynomials}

\author{Fernando De Ter\'an\thanks{Departamento de Matem\'aticas, Universidad Carlos III de Madrid, Avda. Universidad 30, 28911 Legan\'es, Spain. {\tt fteran@math.uc3m.es}} 
\and Froil\'an M. Dopico\thanks{Departamento de Matem\'aticas, Universidad Carlos III de Madrid, Avda. Universidad 30, 28911 Legan\'es, Spain. {\tt dopico@math.uc3m.es}} 
\and Vadym Koval\thanks{Faculty of Mathematics and Computer Science, Jagiellonian University, \L ojasiewicza 6, 30-348, Krak\'ow, Poland. {\tt vadym2312@gmail.com}} 
\and Patryk Pagacz\thanks{Faculty of Mathematics and Computer Science, Jagiellonian University, \L ojasiewicza 6, 30-348, Krak\'ow, Poland.  {\tt patryk.pagacz@gmail.com}}
}

\begin{document}

\date{\today}

\maketitle

\begin{abstract}
Bundles of matrix polynomials are sets of matrix polynomials with the same size and grade and the same eigenstructure up to the specific values of the eigenvalues. It is known that the closure of the bundle of a pencil $L$ (namely, a matrix polynomial of grade $1$), denoted by $\cB(L)$, is the union of $\cB(L)$ itself with a finite number of other bundles. The first main contribution of this paper is to prove that the dimension of each of these bundles is strictly smaller than the dimension of $\cB(L)$. The second main contribution is to prove that also the closure of the bundle of a matrix polynomial of grade larger than $1$ is the union of the bundle itself with a finite number of other bundles of smaller dimension. To get these results we obtain a formula for the (co)dimension of the bundle of a matrix pencil in terms of the Weyr characteristics of the partial multiplicities of the eigenvalues and of the (left and right) minimal indices, and we provide a characterization for the inclusion relationship between the closures of two bundles of matrix polynomials of the same size and grade.   
\end{abstract}



\section{Introduction}

Let $d$ be a non-negative integer. The set of $m\times n$ complex matrix polynomials of grade $d$, or, equivalently, of degree at most $d$, is the set whose elements are of the form
\begin{equation}\label{poly}
    P(\la)=A_0+\la A_1+\cdots+\la^d A_d,
\end{equation}
where $A_0,\hdots,A_d \in \mathbb{C}^{m\times n}$ are complex $m\times n$ matrices, and $A_d$ is allowed to be equal to zero and $\lambda$ is a complex scalar variable. The largest index $k$ such that $A_k\neq0$ in \eqref{poly} is the degree of $P(\la)$. The degree of the identically zero matrix polynomial, i.e., that with $A_0 = A_1 = \cdots = A_d = 0$ is defined to be $-\infty$. Matrix polynomials with grade $d=1$ are called {\em matrix pencils} (or just ``pencils" for short). Observe that the set of $m\times n$ complex matrix polynomials of grade $d$ is a vector space over $\mathbb{C}$.

The {\em eigenstructure} of matrix polynomials plays a relevant role in many of the applied problems where these polynomials arise \cite{MMMM06b,MeVo04,TiMe01}. The subset of the vector space of $m\times n$ complex matrix polynomials of grade $d$ formed by the matrix polynomials having the same eigenstructure form an {\em orbit}, and if the distinct eigenvalues are not fixed, then they form a {\em bundle} instead (see Section \ref{notation_sec} for more details on these notions). The number of different bundles in the vector space of $m\times n$ complex matrix polynomials with grade $d$ is finite, and this space is, then, equal to the finite union of these bundles. The analysis of the inclusion relationships between the closures of orbits and bundles of matrix pencils and matrix polynomials has been a topic of research since, at least, the 1980's, and follows a previous line of research for $n\times n$ matrices that had contributions by eminent mathematicians, like V. I. Arnold \cite{Arno71}. These inclusion relationships allow us to describe the set of $m\times n$ matrix pencils or matrix polynomials with grade $d$ from a theoretical and a topological point of view, which makes the topic to be interesting by itself. Besides, it also has implications in the numerical computation of the eigenstructure of matrix pencils and matrix polynomials, as we are going to see. Actually, this topic is closely connected to the spectral perturbation theory of matrix pencils and matrix polynomials since, if a matrix polynomial, say $Q$, belongs to the closure of the bundle of another matrix polynomial, say $P$ (and $Q$ does not belong to the bundle of $P$), then, in every neighborhood of $Q$, there are matrix polynomials with the same eigenstructure as $P$ (up to the specific values of the distinct eigenvalues), which is different to the one of $Q$. So, an arbitrarily small perturbation of $Q$ can lead to a matrix polynomial with a different eigenstructure, namely that of $P$.

Let us summarize the most relevant contributions so far in the topic of describing the inclusion relationships between orbit and bundle closures of general matrix pencils and matrix polynomials:
\begin{itemize}
    \item Pokrzywa obtained in \cite{Pokr86} a characterization for the inclusion relationships between orbit closures of matrix pencils, that was later reformulated in \cite{Hoyo90} in terms of majorizations of the Weyr characteristics of the eigenstructure of the corresponding pencils (see Section \ref{notation_sec} for more information on these notions, and Theorem \ref{dehoyos_th} for this characterization). A recent new proof of Pokrzywa's result can be found in \cite{sergeichuk2021}.
    \item The most comprehensive work so far regarding the description of the geometry of the set of $m\times n$ matrix pencils from the point of view of orbits and bundles is contained in \cite{EdEK97} and \cite{EdEK99}. In these papers, a characterization of the covering between orbit and bundle closures is presented ({\em covering} refers to an orbit or bundle closure which is included in another one and there is no any other orbit or bundle closure in between). This allows us to obtain the complete Hasse diagram of the inclusion relation between closures of either orbits or bundles of matrix pencils (the so-called {\em complete stratification}). As a consequence of the analysis carried out in these papers, the authors propose, in the second one, an improvement of an algorithm for computing the eigenstructure of $m\times n$ matrix pencils, namely the GUPTRI algorithm introduced in \cite{DeKa93} and \cite{DeKa93b}, that takes into account the stratification of orbit and bundle closures. 
    \item The works \cite{EdEK97} and \cite{EdEK99} led to a series of contributions that aimed to develop a software tool to obtain and display the stratification of sets of matrix pencils. A relevant output produced by this line of research is the Java-based tool {\tt Stratigraph}, that displays the complete stratification of orbits and bundles of $n\times n$ matrices (under similarity) and $m\times n$ matrix pencils (under strict equivalence) for small to moderate values of $m$ and $n$ (see, for instance, \cite{ElJK01}).
    \item The recent work \cite{DJKV20} contains a description of the inclusion relationships of orbit and bundle closures of matrix polynomials with grade larger than $1$ by means of the corresponding ones for pencils using some particular linearizations (namely, matrix pencils associated with matrix polynomials of grade larger than $1$ which allow to easily recover the eigenstructure of the matrix polynomials).  It is an extension of the earlier work \cite{JoKV13}, where the stratification of the set of full rank matrix polynomials was described.
\end{itemize}
We want to emphasize that one of the reasons why the previous references dealt first with matrix pencils is due to the fact that the orbit of a pencil is an orbit in the standard sense (namely, the orbit of an element of a given set under the action of a group), whereas for matrix polynomials this is no longer the case (see Section \ref{notation_sec} for more details).

Nowadays, there is still an active area of research devoted to describing the inclusion relationships between the closures of orbits and bundles of matrix pencils and matrix polynomials which are {\em structured}. By ``structured" we mean those matrix pencils and matrix polynomials whose coefficient matrices enjoy some symmetries which are of interest from the point of view of applications \cite{MMMM06b}, since they typically come from symmetries associated to the physical systems that are modeled by the matrix polynomials in these applications. In particular, the complete stratification of skew-symmetric pencils and skew-symmetric matrix polynomials \eqref{poly} of odd grade larger than $1$ (namely, $A^\top_i=-A_i$, for all $0\leq i\leq d$, where $(\cdot)^\top$ denotes the transpose) have been provided in \cite{DmKa14} and \cite{Dmyt17}, respectively. 

Another current research topic related to this paper is to identify the {\em generic eigenstructures} (namely, those such that the closures of the corresponding bundles contain all the pencils or matrix polynomials with the required properties) of matrix pencils and matrix polynomials with bounded rank, for general and structured matrix pencils and matrix polynomials of grade larger than $1$. More precisely, the generic eigenstructures of matrix pencils with bounded rank were provided in \cite{DeDo08}, and the generic eigenstructures of matrix polynomials with grade larger than $1$ and bounded rank have been given more recently in \cite{DmDo17}. For structured matrix pencils with bounded rank, like symmetric (namely, $A_i^\top=A_i$), Hermitian (namely, $A_i^*=A_i$, where $(\cdot)^*$ denotes the conjugate transpose), $\top$-palindromic (namely, $A_i^\top=A_{d-i}$), $\top$-anti-palindromic (namely $A_i^\top=-A_{d-i}$), $\top$-even (namely $A_i^\top=A_i$ for $i$ even and $A_i^\top=-A_i^\top$ for $i$ odd), and $\top$-odd (namely, $A_i^\top=-A_i$ for $i$ even and $A_i^\top=A_i$ for $i$ odd), the generic eigenstructures have been identified in \cite{ddd20}, \cite{ddd24}, and \cite{DeTe17} respectively. For symmetric matrix polynomials with bounded rank and odd grade, the generic eigenstructures have been also described in \cite{ddd20-2}, while skew-symmetric matrix polynomials with bounded rank are covered in \cite{de2023even,DmDo18} for any grade. However, in this work we only consider general (non-structured and not of bounded rank) matrix pencils and matrix polynomials.

The results about stratifications previously mentioned and published in \cite{Dmyt17,DJKV20,DmKa14,EdEK97,EdEK99} describe the complete inclusion relationships between the closures of orbits and bundles of the corresponding sets of matrix pencils and matrix polynomials. From a geometrical point of view, the notion of ``stratification'' refers to prove that the corresponding set of matrix pencils or matrix polynomials is  ``a stratified manifold''. Quoting \cite[p. 670]{EdEK99}, ``a stratified manifold is the union of nonintersecting manifolds whose
closure is the finite union of itself with strata of {\em smaller dimensions}". In the framework of matrix pencils and matrix polynomials considered in this paper, the strata can either be the orbits or the bundles. In the case of orbits of matrix pencils under strict equivalence, general results \cite[Proposition, p. 60]{Hump75} about orbits of elements under the action of groups guarantee  that those orbits in the closure of a given orbit $\orb (L)$ which are different from $\orb (L)$ have indeed strictly lower dimension than $\orb (L)$. However, for bundles of matrix pencils and matrix polynomials such a result about closures and dimensions is not explicitly proved in the literature. To fill this gap is the main goal of this paper.

Let $P$ be an $m\times n$ matrix polynomial of grade $d$ and let $\cB(P)$ be the bundle of $P$, namely the set of  $m\times n$ matrix polynomials of grade $d$ with the same eigenstructure as $P$, up to the specific values of the distinct eigenvalues. It has been recently proved in \cite[Theorem 16]{DD23} that for $d=1$, i.e., for matrix pencils, the closure of $\cB(L)$ is the finite union of $\cB(L)$ itself together with other bundles, but it remained as an open problem to prove that the dimension of these other bundles is strictly smaller than the dimension of $\cB(L)$ (see \cite[Section 5]{DD23}). In addition, for matrix polynomials $P$ of grade larger than $1$, it remains as an open problem to prove that the closure of $\cB(P)$ is the finite union of $\cB(P)$ itself together with other bundles, and also that these other bundles have dimensions strictly smaller than $\cB(P)$.

The main contributions of this paper, within the previous framework, are the following:
\begin{itemize}
    \item To prove that for any pencil $L$ the closure of $\cB(L)$ is the finite union of $\cB(L)$ together with other bundles with strictly smaller dimension (Theorem \ref{T16}). This was mentioned as an open problem in \cite[Section 5]{DD23}.
    \item To prove that for any matrix polynomial $P$ of grade larger than $1$ the closure of $\cB(P)$ is the finite union of $\cB(P)$ together with other bundles with strictly smaller dimension (Theorem \ref{union_th}). 
\end{itemize}
We also provide a new formula for the codimension of the orbit of a matrix pencil $L$ in terms of the Weyr characteristics of the eigenstructure of $L$ (Proposition \ref{orbitcodim_prop}), which is key to prove Theorem \ref{T16}. We also provide a characterization for a matrix polynomial of grade larger than $1$ to belong to the closure of the bundle of another matrix polynomial of the same size and grade (Theorem \ref{main_th}), which extends the characterization for pencils provided in \cite[Th. 12]{DD23} and solves another open problem posed in \cite[Section 5]{DD23}.

The rest of the paper is organized as follows. In Section \ref{notation_sec} we recall the notation and basic notions, together with some elementary results, that are used throughout the paper. In Section \ref{weyr_sec} we present the formula for the codimension of the orbit of a matrix pencil $L$ in terms of the Weyr characteristics of the eigenstructure of $L$. Section \ref{main_sec} is devoted to present the main results of this paper. More precisely, in Section \ref{pencils_sec} we prove that the dimension of any bundle in the closure of $\cB(L)$, other than $\cB(L)$, is strictly smaller than the dimension of $\cB(L)$, whereas Section \ref{closures_sec} is devoted to prove that the closure of $\cB(P)$, with $P$ being a matrix polynomial of grade larger than $1$, is the finite union of $\cB(P)$ together with other bundles of strictly smaller dimension. Finally, in Section \ref{conclusions_sec} we summarize the main contributions of the paper and discuss some lines of possible future research related with these contributions.

\section{Basic notions and notation}\label{notation_sec}

We use the following notation throughout the paper: $I_n$ stands for the $n\times n$ identity matrix, $\CC^{m\times n}$ stands for the set of $m\times n$ matrices with complex entries, and  $\overline\CC=\CC\cup\{\infty\}$. 

For brevity, and when there is no risk of confusion, we often omit the reference to the variable $\la$ in the matrix polynomial $P(\la)$ and we just write $P$. For matrix polynomials of arbitrary grade we use the letters $P$ and $Q$, whereas for matrix pencils we use $L$ and $M$ instead. All the matrix pencils and matrix polynomials considered in this paper have complex coefficients. The normal rank, or simply the rank, of a matrix polynomial $P(\lambda)$ is denoted as $\rank P$ and is the size of the largest non-identically zero minor of $P(\lambda)$. A matrix polynomial $P(\lambda)$ is regular if is square and $\det P(\lambda) \not\equiv 0$. Otherwise $P(\lambda)$ is singular. 

Two matrix polynomials $P(\la)$ and $Q(\la)$ are said to be {\em strictly equivalent} if there are two (constant) invertible matrices $R$ and $S$ such that $Q(\la)=R P(\la)S$. 

The {\em eigenstructure} of a matrix polynomial $P$ of grade $d$ consists of the set of {\em left} and {\em right minimal indices}, together with the {\em partial multiplicities} associated to the finite and infinite eigenvalues (for the definition of these notions see, for instance, \cite{DeDM14} and recall that the definition of the infinite eigenvalue depends on the grade $d$). We use the notation $\Lambda(P)$ for the spectrum of $P$ (namely, the set of distinct eigenvalues, possibly including the infinite one).

The eigenstructure and the strict equivalence are related notions for matrix pencils. More precisely, two matrix pencils $L$ and $M$ are strictly equivalent if and only if they have the same eigenstructure, which is encoded in the so-called {\em Kronecker Canonical Form} (KCF) \cite[Ch. XII, \S5, Th. 5]{Gant59}. In particular, the partial multiplicities are associated to {\em Jordan blocks} corresponding to eigenvalues in the KCF (there is a block of size $s\times s$ with eigenvalue $\la$ for each partial multiplicity $s$ corresponding to $\la$), whereas right (respectively, left) minimal indices are associated to {\em right} (respectively, {\em left}) {\em singular blocks} with size $r\times(r+1)$ (respectively, $(\ell+1)\times\ell$), where $r$ (resp., $\ell$) is the corresponding right (resp., left) minimal index.

The left and the right minimal indices of an $m\times n$ matrix polynomial $P$ form lists on nonnegative integers, while, for each eigenvalue of $P$, its partial multiplicities form a list of positive integers, and these lists are known as the {\em Segre characteristics}. If $\mu\in\overline{\CC}$ is not an eigenvalue of $P$, we say that $P$ has not partial multiplicities at $\mu$. Instead of working with the Segre characteristics, we deal with the so-called {\em Weyr characteristics}, which are the conjugate lists (see, for instance, \cite{weyr_book} and \cite{shapiro} for a more comprehensive analysis of the Weyr characteristic, and its relation with the Segre characteristic, in the case of matrices instead of matrix pencils). More precisely:
\begin{itemize}
    \item $r(P):=(r_0(P),r_1(P),\hdots)$ is the Weyr characteristic of right minimal indices of $P$, where $r_i(P)$ is the number of right minimal indices of $P$ which are greater than or equal to $i$;
    \item $\ell(P):=(\ell_0(P),\ell_1(P),\hdots)$ is the Weyr characteristic of the left minimal indices of $P$, where $\ell_i(P)$ is the number of left minimal indices of $P$ which are greater than or equal to $i$; and
    \item $W(\mu,P):=(W_1(\mu,P),W_2(\mu,P),\hdots)$ is the Weyr characteristic of the partial multiplicities of $P$ at $\mu\in\overline{\CC}$, where $W_i(\mu,P)$ is the number of partial multiplicities of $P$ associated with $\mu$ which are greater than or equal to $i$. Observe that the terms of $W(\mu , P)$ are {\em all} equal to zero if and only if $\mu$ is not an eigenvalue of $P$.
\end{itemize}
Note that all three lists $r(P),\,\ell(P)$, and $W(\mu,P)$ have a finite number of nonzero elements. We can, therefore, consider them either as finite lists, after removing all zero terms, or as infinite lists all whose terms are zero from a certain index on. We will use the following two equalities $\rank P + r_0 (P) = n$ and $\rank P + \ell_0 (P) = m$, which follow from the rank-nullity theorem. Observe that a matrix polynomial is regular if and only if it has no minimal indices, neither left nor right.

Some majorizations of the previous Weyr characteristics are key in a relevant part of this work. We recall that the list $(a_1,a_2,\hdots)$, where $a_1 \geq a_2 \geq \cdots$, {\em majorizes} the list $(b_1,b_2,\hdots)$, where $b_1 \geq b_2 \geq \cdots$, denoted by $(b_1,b_2,\hdots)\prec(a_1,a_2,\hdots)$, if $\sum_{i=1}^jb_i\leq\sum_{i=1}^ja_i$, for all $j\geq1$.  

We need the following notions, that have been used in previous references (see, for instance, \cite{DD23,DeDM12,DeDM14,DmDo17}).
\begin{itemize}
    \item $\cC_P$ denotes the {\em first Frobenius companion linearization} of the matrix polynomial of grade $d$, $P(\la)=\sum_{i=0}^d\la^{i}A_i$, with $A_i\in\CC^{m\times n}$, for $0\leq i\leq d$ (see, for instance, \cite[\S5.1]{DeDM14}), that is
    $$
    \cC_P(\la):=\la \begin{bmatrix}
        A_d&&&\\&I_n&&\\&&\ddots&\\&&&I_n
    \end{bmatrix}+\begin{bmatrix}
        A_{d-1}&A_{d-2}&\cdots&A_0\\-I_n&0&\cdots&0\\&\ddots&\ddots&\vdots\\0&&-I_n&0
    \end{bmatrix}.
    $$
    We will refer to these pencils as ``companion pencils" for short.
    
    It is worth to emphasize that there is a one-to-one correspondence between the eigenstructure of $P$, considered as a matrix polynomial of grade $d$, and that of $\cC_P$, considered as a matrix pencil, namely, as a polynomial of grade $1$ (see, for instance, \cite[Th. 5.3]{DeDM14}). More precisely:
    \begin{itemize}
        \item Both $P$ and $\cC_P$ have the same eigenvalues (finite and infinite), with the same partial multiplicities.
        \item If $\varepsilon_1\leq\cdots\leq\varepsilon_p$ are the right minimal indices of $P$, then $\varepsilon_1+d-1\leq\cdots\leq\varepsilon_p+d-1$ are the right minimal indices of $\cC_P$.
        \item If $\eta_1\leq\cdots\leq\eta_q$ are the left minimal indices of $P$, then $\eta_1\leq\cdots\leq\eta_q$ are the left minimal indices of $\cC_P$.
    \end{itemize}    
   \item ${\rm POL}_d^{m\times n}$ denotes the set of matrix polynomials of size $m\times n$ and grade $d$.
    \item ${\rm GSYL}_d^{m\times n}$ denotes the set of matrix pencils of the form ${\cal C}_P$, with $P\in{\rm POL}_d^{m\times n}$ (this is the so-called {\em generalized Sylvester space} in \cite{DmDo17}).
     \item The {\it orbit} of $P\in{\rm POL}_d^{m\times n}$, denoted by $\orb(P)$, is the set of $m\times n$ matrix polynomials of grade $d$ having exactly the same eigenstructure as $P$. %All elements of $\orb(P)$ have the same degree as $P$.
     In the case of a pencil $L$, its orbit, $\orb(L)$, is the set of matrix pencils which are strictly equivalent to $L$, so, in this case, it coincides with the standard notion of orbit of an element of a set (the set of matrix pencils) under the action of a group (namely, ${\rm GL}_m(\CC)\times{\rm GL}_n(\CC)$, where ${\rm GL}_n(\CC)$ denotes the group of $n\times n$ invertible matrices with complex entries), that is 
     $$
     \begin{array}{ccc}
          {\rm GL}_m(\CC)\times{\rm GL}_n(\CC)\times{\rm POL}_1^{m\times n}& \rightarrow&{\rm POL}_1^{m\times n}\\
          (P,Q,A+\la B)&\mapsto&P(A+\la B)Q 
     \end{array}.
     $$
     For matrix polynomials of grade larger than $1$, this is not the case anymore, since there are matrix polynomials of the same size and grade having the same eigenstructure but not being strictly equivalent (see, for instance, \cite[p. 277]{DeDM14}).
    \item For $P \in {\rm POL}_d^{m\times n}$, the {\em Sylvester orbit of $\cC_P$} is defined in \cite{DmDo17} as $\orb^{\rm syl}(\cC_P):=\orb(\cC_P)\cap{\rm GSYL}_d^{m\times n}$. Namely, $\orb^{\rm syl}(\cC_P)$ contains all pencils which are strictly equivalent to $\cC_P$ and which are companion pencils of some $m\times n$ matrix polynomial of grade $d$. 
    \item The {\it bundle} of $P\in{\rm POL}_d^{m\times n}$, denoted by $\cB(P)$, is the set of $m\times n$ matrix polynomials of grade $d$ having exactly the same eigenstructure as $P$, up to the specific values of the distinct eigenvalues, i.e., as long as those eigenvalues which are distinct in $P$ remain distinct in $\widetilde P\in\cB(P)$.    
    \item For $P \in {\rm POL}_d^{m\times n}$, the {\em Sylvester bundle of $\cC_P$} is defined in \cite{DD23} as $\cB^{\rm syl}(\cC_P):=\cB(\cC_P)\cap{\rm GSYL}_d^{m\times n}$. Namely, $\cB^{\rm syl}(\cC_P)$ is the set of matrix pencils in $\cB(\cC_P)$ which are companion pencils for some $m\times n$ matrix polynomial of grade $d$.
\end{itemize}
We also  
note the following definition, which extends Definition 1 in \cite{DD23} from matrix pencils to matrix polynomials.

\begin{definition}\label{coalescence_def}
Let $P  \in {\rm POL}_d^{m\times n}$ with distinct eigenvalues $\mu_1,\hdots,\mu_s$, and let $\psi:\overline\CC\rightarrow\overline\CC$ be a map. Then, $\psi_c(P)  \in {\rm POL}_d^{m\times n}$ is any matrix polynomial satisfying the following three properties:
\begin{itemize}
    \item $r(\psi_c(P))=r(P)$,
    \item $\ell(\psi_c(P))=\ell(P)$, and
    \item $W(\mu,\psi_c(P))=\bigcup_{\mu_i\in\psi^{-1}(\mu)}W(\mu_i,P)$, for all $\mu\in\overline{\CC}$.
\end{itemize}
\end{definition}

In {\rm\cite{DD23}} it is said that a pencil $\psi_c(L)$, as in Definition {\rm\ref{coalescence_def}}, is obtained from $L$ after {\em coalescing} some eigenvalues of $L$. Recall, see \cite{DD23}, that the union of Weyr characteristics above is carried out with repetition of entries.

As it is emphasized in \cite[Remark 2]{DD23}, if $P \ne 0$, there are infinitely many matrix polynomials $\psi_c(P)$ as in Definition \ref{coalescence_def}, but all of them have the same eigenstructure and, so, they are all in the same orbit. Moreover, note that the equality
\begin{equation}\label{psi-comp}
\orb (\psi_c(\cC_P))=\orb (\cC_{\psi_c(P)})
\end{equation}
follows immediately from the relation between the eigenstructures of $\psi_c(P)$ and $\cC_{\psi_c(P)}$.

If $P\in{\rm POL}_d^{m\times n}$, then $\overline\orb(P)$ and $\overline\cB(P)$ denote the closure of, respectively, the orbit and the bundle of $P$ in the Euclidean topology of ${\rm POL}_d^{m\times n}$ identified with $\CC^{(d+1)mn}$.

Let us recall the following map introduced in \cite{DmDo17},
\begin{equation}\label{homeo}
\begin{array}{cccc}
     f:&{\rm POL}_d^{m\times n}&\rightarrow&{\rm GSYL}_d^{m\times n}  \\
     &P&\mapsto&{\cal C}_P,
\end{array}
\end{equation}
which is a homeomorphism (in fact, a bijective isometry, see \cite[p. 218]{DmDo17}) satisfying 
\begin{equation}\label{forb}
f^{-1}(\orb^{\rm syl}(\cC_P))=\orb(P)
\end{equation} 
(or, equivalently, $f(\orb(P))=\orb^{\rm syl}(\cC_P)$), 
 for any $P\in {\rm POL}_d^{m\times n}$, which is equivalent to say that $\orb^{\rm syl}(\cC_P)$ is the set of all companion pencils of matrix polynomials in ${\rm POL}_d^{m\times n}$ having the same eigenstructure as $P$. This follows immediately from the one-to-one correspondence between the eigenstructures of $P$ and $\cC_P$ mentioned before. Moreover, since $f$ is a homeomorphism,
\begin{equation}\label{invforb}
    f^{-1}(\overline{\orb^{\rm syl}}(\cC_P))=\overline\orb(P)
\end{equation}
(see \cite[Lemma 2.4]{DmDo17}). The identities \eqref{forb} and \eqref{invforb} can be extended to bundles. More precisely, $f^{-1}(\cB^{\rm syl}(\cC_P))=\cB(P)$, which is, again, an immediate consequence of the one-to-one correspondence between the eigenstructures of $P$ and $\cC_P$, and this identity implies, because $f$ is a homeomorphism, that
\begin{equation}\label{bunclosuref}
    f^{-1}(\overline{\cB^{\rm syl}}(\cC_P))=\overline\cB(P)
\end{equation}
(see also \cite[p. 1047]{ddd20-2}, where the corresponding identities for symmetric matrix polynomials of odd grade and symmetric companion forms are presented). The identity \eqref{bunclosuref} immediately implies that
\begin{equation}\label{bunequivalence}
\overline{\cB}(P)\subseteq\overline\cB(Q)\qquad\mbox{if and only if}\qquad \overline{\cB^{\rm syl}}(\cC_P)\subseteq\overline{\cB^{\rm syl}}(\cC_Q)
\end{equation}
(the corresponding equivalence for symmetric matrix polynomials of odd grade and symmetric companion forms is also presented in \cite[p. 1047]{ddd20-2}).

The following lemma is the counterpart of Lemma 3.1 in \cite{DmDo17} for bundles instead of orbits, and also the one of Lemma 6.2 in \cite{ddd20-2} for unstructured matrix polynomials instead of symmetric ones.

\begin{lemma}\label{sylvclosure_lem}
    Let $P(\la)$ be an $m\times n$ matrix polynomial with grade $d$. Then
    $$
    \overline{\cB^{\rm syl}}(\cC_P)=\overline{\cB}(\cC_P)\cap{\rm GSYL}_d^{m\times n}.
    $$
\end{lemma}
\begin{proof} Since the closure of the intersection of two sets is always included in the intersection of the closures of both sets, the inclusion $\overline{\cB^{\rm syl}}(\cC_P)\subseteq\overline{\cB}(\cC_P)\cap{\rm GSYL}_d^{m\times n}$ follows (take into account that ${\rm GSYL}_d^{m\times n}$ is closed). 

Let us now prove the reverse inclusion, so let $L\in \overline{\cB}(\cC_P)\cap{\rm GSYL}_d^{m\times n}$. As $L\in{\rm GSYL}_d^{m\times n}$, it must be $L=\cC_{\widetilde P}$, for some $m\times n$ matrix polynomial $\widetilde P$ of grade $d$. Since $L\in\overline\cB(\cC_P)$, there is a sequence of pencils $\{L_i\}_{i\in\mathbb N}$, with $L_i\in\cB(\cC_P)$, for all $i\in\mathbb N$, that converges to $L=\cC_{\widetilde P}$. By Theorem 2.5 in \cite{DmDo17}, for $L_i$ close enough to $\cC_{\widetilde P}$, the pencil $L_i$ is strictly equivalent to a companion pencil $\cC_{P_i}$ which is also close to $\cC_{\widetilde P}$, namely the sequence $\{\cC_{P_i}\}_{i\in\mathbb N}$ converges to $\cC_{\widetilde P}$. Since $\cC_{P_i}\in\cB(\cC_P)$, we conclude that  $L=\cC_{\widetilde P}\in\overline{ \cB(\cC_P)\cap{\rm GSYL}_d^{m\times n}}$, and we are done. 
\end{proof}

Also, for $f$ being the map in \eqref{homeo}, 
    \begin{equation}\label{fclosure}
    f(\overline{\cB}(P))=\overline{\cB^{\rm syl}}(\cC_P)=%\overline{\cB(\cC_P)\cap{\rm GSYL}_d^{m\times n}}=
    \overline\cB(\cC_P)\cap{\rm GSYL}_d^{m\times n},
    \end{equation}
 where the first identity is equivalent to \eqref{bunclosuref}, and the second one is Lemma \ref{sylvclosure_lem}.

The following lemma is the generalization of Lemma 11 in \cite{DD23} from matrix pencils to matrix polynomials of grade larger than $1$.

\begin{lemma}\label{pinclosureq_lemma}
    Let $P,Q\in{\rm POL}_d^{m\times n}$. Then $P\in\overline\cB(Q)$ if and only if $\overline\cB(P)\subseteq\overline\cB(Q)$. 
\end{lemma}
\begin{proof}
The inclusion $\overline\cB(P)\subseteq\overline\cB(Q)$ immediately implies $P\in\overline\cB(Q)$, since $P\in\cB(P)\subseteq \overline\cB(P)$.

If $P\in\overline\cB(Q)$, from \eqref{fclosure} we get $\cC_P=f(P)\in f(\overline{\cB}(Q))=\overline{\cB^{\rm syl}}(\cC_Q)=\overline\cB(\cC_Q)\cap{\rm GSYL}_d^{m\times n}$. In particular, $\cC_P\in\overline\cB(\cC_Q)$, so \cite[Lemma 11]{DD23} implies that $\overline\cB(\cC_P)\subseteq\overline\cB(\cC_Q)$, and this in turn implies $\overline{\cB^{\rm syl}}(\cC_P)\subseteq\overline{\cB^{\rm syl}}(\cC_Q)$, by Lemma \ref{sylvclosure_lem}. Now, by \eqref{bunclosuref}, we get $\overline{\cB}(P)=f^{-1}(\overline{\cB^{\rm syl}}(\cC_P))\subseteq f^{-1}(\overline{\cB^{\rm syl}}(\cC_Q))=\overline\cB(Q)$, and this concludes the proof.
\end{proof}

It is known that $\orb(L)$, with $L$ being an $m\times n$ matrix pencil, is a differentiable manifold (see, for instance, \cite{DeEd95}), so its {\em dimension} is the dimension as a differentiable manifold (namely, the dimension of the tangent space). Its {\em codimension} is the codimension when considered as a manifold in the space of $m\times n$ matrix pencils, namely ${\rm codim}\,\orb(L)=2mn-\dim\orb(L)$. For a matrix polynomial $P\in{\rm POL}_d^{m\times n}$, instead, the codimension of $\orb(P)$ is defined as the codimension of $\orb^{\rm syl}(\cC_P)$ in the set ${\rm GSYL}_d^{m\times n}$.  It is known \cite[Lemma 1]{DJKV20} that
\begin{equation}\label{codimpoly}
    {\rm codim}\, \orb^{\rm syl}(\cC_P)={\rm codim}\,\orb(\cC_P)
\end{equation}
(note that the first codimension in \eqref{codimpoly} refers to the codimension in the space ${\rm GSYL}_d^{m\times n}$ of all companion pencils of $m\times n$ matrix polynomials of grade $d$, whereas the second one refers to the codimension in the whole space of $m\times n$ matrix pencils).

As for bundles, we follow the definition of codimension provided in previous references, like, for instance, \cite[p. 443]{DJKV20}, namely
\begin{equation}\label{codimbundle_def}
    {\rm codim}\,\cB(P)={\rm codim}\,\orb(P)-\#\Lambda(P),
\end{equation}
for any $P\in{\rm POL}_d^{m\times n}$. Therefore

\begin{equation}\label{codimbcp}
\begin{array}{ccl}
    \codim\cB^{\rm syl}(\cC_P)&=&\codim\orb^{\rm syl}(\cC_P)-\#\Lambda(\cC_P)\\
    &=&\codim\orb(P)-\#\Lambda(P)\\&=&\codim\cB(P).
    \end{array}
\end{equation}

\section{Formula for the codimension of the orbit in terms of Weyr characteristics}\label{weyr_sec}

A formula for the codimension of the orbit of a matrix pencil $L$ is provided in \cite{DeEd95} in terms of the minimal indices and the partial multiplicities of $L$. Here we provide an alternative formula using the Weyr characteristics.

\begin{proposition}\label{orbitcodim_prop}
Let $L(\lambda)$ be an $m\times n$ matrix pencil with distinct eigenvalues $\lambda_1,\hdots,\la_p\, \in \overline\CC$. Then 
\begin{equation}\label{codimorbit}
\begin{array}{ccl}
\textnormal{codim } \mathcal{O}(L)&=&\displaystyle \ell_0(L)n+r_0(L)m-\displaystyle\sum_{i=0}^\infty{r_i(L)r_{i+1}(L)}\\&&\displaystyle-\sum_{i=0}^\infty{\ell_i(L)\ell_{i+1}(L)}+\sum_{k=1}^p\sum_{i=1}^\infty{W_i(\lambda_k,L)^2}.
\end{array}
\end{equation}
\end{proposition}
\begin{proof}
The calculation of the codimension of the orbit of a pencil in \cite{DeEd95} is divided in 5 parts. We will translate each of them separately in terms of the Weyr characteristics. For brevity, we remove in this proof the reference to the pencil $L$ in the Weyr characteristics, and write $r_i,\ell_i$, and $W_i(\mu)$ instead of $r_i(L),\ell_i(L)$, and $W_i(\mu,L)$, respectively.

\medskip

1. {\em Codimension of the Jordan structure}. 

Let $q_1(\mu) \geq q_2(\mu)\geq  q_3(\mu)\geq\cdots$ denote the sizes of the Jordan blocks in the KCF of
$L$ corresponding to $\mu\in\overline{\CC}$. Then, the Segre characteristic of $\mu$ is the list $q(\mu) = (q_1 (\mu) , q_2 (\mu) , q_3 (\mu), \ldots )$. The
codimension of the Jordan structure of $L$ is expressed in \cite{DeEd95} as 
\begin{equation}\label{jordansum}
    \sum_{\mu\in\Lambda(L)} \sum_{i=1}^\infty (2i-1)q_i(\mu).
\end{equation}

For each $\mu$ let us consider the Ferrers diagram, i.e. a set of weighted cells organized in columns, where the $i$th column has height $q_i(\mu)$ and the weight of each cell in this column is equal to $2i-1$ (see Figure \ref{Fig}).  Then the sum of all weights in a given column is equal to $(2i-1)q_i(\mu)$. Therefore, the sum of all weights in the diagram is 
$
q_1(\mu)+3q_2(\mu)+5q_3(\mu)+\cdots
$ 

\begin{figure}
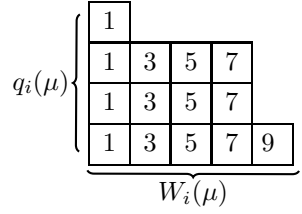

    \centering
    
\begin{ytableau}
    \tikznode{a3} 1 \\
      1 & 3 & 5 & 7 \\
      1 & 3 & 5 & 7 \\
    \tikznode{a1} 1 & 3 & 5 & 7 & 9  \tikznode{a2}{~} & \none & \none 
\end{ytableau}

\tikz[overlay,remember picture]{%
\draw[decorate,decoration={brace},thick] ([yshift=-2.4mm,xshift=2mm]a2.south east) -- 
([yshift=-2.5mm,xshift=-2mm]a1.south west) node[midway,below]{$W_i(\mu)$};
}
\tikz[overlay,remember picture]{%
\draw[decorate,decoration={brace},thick] ([yshift=0mm,xshift=-5mm]a1.south east) -- 
([yshift=2mm,xshift=-3mm]a3.south west) node[midway,left]{$q_i(\mu)$};
}
\bigskip
\caption{Ferrers diagram for an eigenvalue $\mu$ with partial multiplicities $4$, $3$, $3$, $3$, and $1$, so $q(\mu)=(4,3,3,3,1)$ and $W(\mu)=(5,4,4,1)$.}\label{Fig}
\end{figure}

Now, note that the length of the $i$th row is $W_i(\mu)$.
Thus, if we compute the sum of all weights in the diagram row-wise (namely, adding up first the weights in all cells in each row and then adding up the resulting numbers for all rows), we get 
$
W_1(\mu)^2+W_2(\mu)^2+\cdots $
Therefore, the total codimension of the Jordan structure \eqref{jordansum} is equal to
\begin{equation}\label{W2}
\sum_{k=1}^p\sum_{i=1}^\infty{W_i(\lambda_k)^2}.
\end{equation}

2. {\em Codimension of the right singular blocks.}

Let $a_1 \geq a_2\geq\cdots$ denote the right minimal indices of $L$. The Segre characteristic of the right minimal indices of $L$ is the list $a = (a_1 , a_2 , \ldots )$. The codimension corresponding to the right minimal indices was expressed in \cite{DeEd95} as 
\begin{equation}\label{ais}
\sum_{a_j>a_k} (a_j-a_k-1).
\end{equation}
We consider a similar Ferrers diagram as in the previous case but with the height of the $i$th column equal to $a_i + 1$ for taking into account that some minimal indices may be equal to zero. In this case, no weights are assigned to the cells of the initial Ferrers diagram (see the leftmost diagram in Figure \ref{Fig2}). Observe that the number of cells of the $i$th row (counting upwards and starting with the zero row at the bottom) is $r_i$, where $r = (r_0, r_1, \ldots )$ is the Weyr characteristic of the right minimal indices of $L$. Next, each $i$th row of this initial diagram except the top one is completed with $r_0 - r_i$ ``fake'' cells. 

Our goal is to fill these ``fake'' cells with weights that sum up the quantity in \eqref{ais}. We emphasize that no weights will be assigned to the cells of the original Ferrers diagram in this process. For this purpose, we initially assign a zero weight to each ``fake'' cell. Then, we start with the index $a_1$ and consider all the other indices $a_k$ such that $a_1 - a_k -1 > 0$. For each of such index $a_k$, we add $+1$ to the weights of the ``fake'' cells in the $k$th column located in rows $a_k + 1, a_k + 2 , \ldots , a_1-1$. At the end of this first step, the sum of the weights of all the ``fake'' cells is $\sum_{a_1>a_k} (a_1-a_k-1)$ (see the second diagram in the first row of the example of Figure \ref{Fig2}).

Next, we consider the index $a_2$ and all the other indices $a_k$ such that $a_2 - a_k -1 > 0$ and, for each of them, we repeat the process of adding $+1$ to the weights of the ``fake'' cells in the $k$th column located in rows $a_k + 1, a_k + 2 , \ldots , a_2-1$. At the end of this second step the sum of the weights of all the ``fake'' cells is $\sum_{a_1>a_k} (a_1-a_k-1) + \sum_{a_2>a_k} (a_2-a_k-1)$ (see the third diagram in the first row of the example of Figure \ref{Fig2}). We repeat this process with all the remaining indices $a_3, a_4, \ldots$, except with the last one. At the end of the process, the sum of all the weights of the ``fake'' cells is precisely \eqref{ais} (see the diagram in the second row of the example of Figure \ref{Fig2}). 

\begin{tikzpicture}[overlay,remember picture]

\draw[red, thick] ([yshift=3mm,xshift=0mm] pic cs:pp) -- ([yshift=-1mm,xshift=0mm] pic cs:pp2);
\end{tikzpicture}

\begin{tikzpicture}[overlay,remember picture]
\draw[red, thick] ([yshift=3mm,xshift=0mm] pic cs:pp3) -- ([yshift=-1mm,xshift=0mm] pic cs:pp4);
\end{tikzpicture}
\begin{tikzpicture}[overlay,remember picture]
\draw[red, thick] ([yshift=3mm,xshift=0mm] pic cs:p) -- ([yshift=-1mm,xshift=0mm] pic cs:p2);
\end{tikzpicture}

\begin{tikzpicture}[overlay,remember picture]
\draw[red, thick] ([yshift=3mm,xshift=2mm] pic cs:p3) -- ([yshift=-1mm,xshift=1mm] pic cs:p4);
\end{tikzpicture}

\begin{figure}
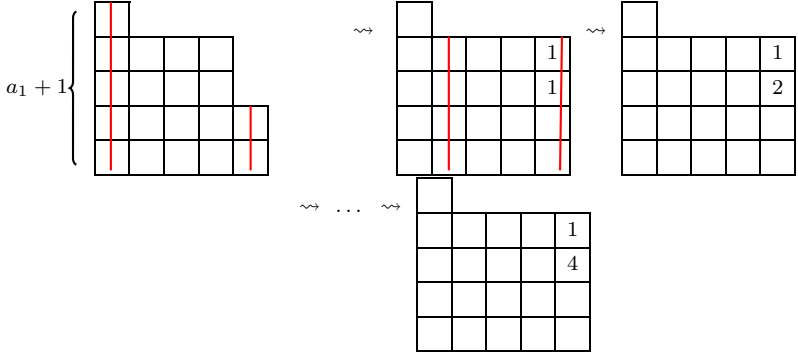

    \centering
   \footnotesize
\begin{ytableau}
    \tikzmark{pp}\tikznode{a3}   \empty \\
       &  &  &  \\
       &  &  &  &   \none \\
   & & & & \tikzmark{pp3} \\
   \tikzmark{pp2} \tikznode{a1}
   \empty &  &  &  &
   \tikzmark{pp4}
    \empty  & \none & \none 
   
\end{ytableau}
~$\leadsto$~
\tikz[overlay,remember picture]{%
\draw[decorate,decoration={brace},thick] ([yshift=0mm,xshift=-5mm]a1.south east) -- 
([yshift=2mm,xshift=-4.5mm]a3.south west) node[midway,left]{$a_1 + 1$};
}
\begin{ytableau}
      \empty \\
     & \tikzmark{p} &  & &  \tikzmark{p3}1 \\ 
       &  &  & & 1
       \\   & & & & \\
     \empty  & \tikzmark{p2} &  &  & \empty \tikzmark{p4}  
\end{ytableau}
~$\leadsto$~
\begin{ytableau}
    \tikznode{a3}   \empty \\
       &  &  & & 1 \\
       &  &  & & 2 \\   & & & & \\
    \tikznode{a1}  \empty &  &  &  & \empty  \tikznode{a2}{~} & \none & \none 
\end{ytableau}
~$\leadsto$~~\dots~~$\leadsto$~
\begin{ytableau}
    \tikznode{a3}   \empty \\
       &  &  & & 1  \\
       &  &  & & 4  \\   & & & & \\
    \tikznode{a1}  \empty &  &  &  & \empty  \tikznode{a2}{~} & \none & \none 
\end{ytableau}
%\bigskip
\caption{Ferrers diagram for right minimal indices $4$, $3$, $3$, $3$, and $1$, including the zero row, so $a=(4,3,3,3,1)$ and $r=(5,5,4,4,1)$ and weights obtained in the last column.}\label{Fig2}
\end{figure}

Moreover, observe that, when the whole process is over, each cell of the $r_0 - r_i$ ``fake'' cells in the $i$th row of the diagram has received a weight $+1$ from each minimal index strictly larger than $i$ and, thus, it has a weight equal to $r_{i+1}$. Therefore, if the weights of the ``fake'' cells are summed up by rows, we get that \eqref{ais} is equal to  $\sum_{i=1}^\infty{(r_0-r_i)r_{i+1}},$ which can be rewritten as
\begin{equation}\label{codimL}
r_0\sum_{i=0}^\infty{r_i}-\sum_{i=0}^\infty{r_ir_{i+1}}-r_0^2.
\end{equation}

\medskip

3. {\em Codimension of the left singular blocks.}

\noindent The codimension due to the left minimal indices provided in \cite{DeEd95} is the same as \eqref{ais} if $a_1\geq a_2\geq\cdots$ correspond to the left minimal indices instead. Proceeding in the same way as for the right minimal indices, we get that the codimension corresponding to this part of the eigenstructure is equal to
\begin{equation}\label{codimLT}
\ell_0\sum_{i=0}^\infty \ell_i-\sum_{i=0}^\infty{\ell_i\ell_{i+1}}-\ell_0^2.    
\end{equation}

\medskip

4. {\em Codimension due to interaction of the Jordan structure with the singular blocks.}

\noindent The codimension of this part provided in \cite{DeEd95} is
$$
(\text{size of the Jordan structure})\cdot (\text{number of singular blocks}),
$$
which can be rewritten as 
\begin{equation}
\label{codimJ}
(r_0+\ell_0)\cdot\sum_{k=1}^p\sum_{i=1}^\infty{W_i(\lambda_k)}.
\end{equation} 

\medskip

5. {\em Codimension due to interactions between left and right singular
blocks.}

\noindent %Let $b_1 \geq b_2\geq b_3\geq\dots$ denote the left
Let $a_1 \geq a_2\geq a_3\geq\cdots$ and $b_1 \geq b_2\geq b_3\geq\cdots$ denote, respectively, the right and left minimal indices of $L$. The term of the codimension corresponding to the interactions between left and right minimal indices obtained in \cite{DeEd95} is equal to
$
\sum_{a_j,b_k}(a_j+b_k+2),
$
that can be rewritten as
\begin{equation}
\label{codimLLT}
\begin{array}{ccl}
     \displaystyle\sum_{b_k}\sum_{a_j}({a_j+b_k+2})&=&\displaystyle\sum_{b_k}\left(b_kr_0+2r_0+{\sum_{i=1}^\infty{r_i}}\right)\\
&=&\displaystyle\sum_{b_k}\left((b_k+1)r_0+{\sum_{i=0}^\infty{r_i}}\right)=\ell_0\sum_{i=0}^\infty{r_i}+r_0\sum_{i=0}^\infty{\ell_i}, 
\end{array}  
\end{equation}
where, to get the first identity, we have used that $\sum_{i=1}^\infty a_i=\sum_{i=1}^\infty r_i$ and, to get the last one, we use that $\sum_{i=1}^\infty b_i=\sum_{i=1}^\infty \ell_i$.

Summing together the expressions from \eqref{codimL}, \eqref{codimLT}, \eqref{codimJ}, and \eqref{codimLLT} we get 
\begin{equation}\label{codimLLT-2}
\begin{array}{c}
\displaystyle(r_0+\ell_0)\left(\sum_{i=0}^\infty{r_i}+\sum_{i=0}^\infty{\ell_i}+\sum_k\sum_{i=1}^\infty{W_i(\lambda_k)}\right)-r_0^2-\ell_0^2-\sum_{i=0}^\infty{r_ir_{i+1}}-\sum_{i=0}^\infty{\ell_i\ell_{i+1}}
\\\displaystyle=(r_0+\ell_0)(n+\ell_0)-r_0^2-\ell_0^2-\sum_{i=0}^\infty{r_ir_{i+1}}-\sum_{i=0}^\infty{\ell_i\ell_{i+1}}.
\end{array}
\end{equation}

Note that
$$
(r_0+\ell_0)(n+\ell_0)-r_0^2-\ell_0^2=n(r_0+\ell_0)+r_0(\ell_0-r_0)=n(r_0+\ell_0)+r_0(m-n)=n\ell_0+mr_0,
$$
since $\ell_0-r_0=m-n$, because the rank of $L$ is equal to both $m-\ell_0$ and $n-r_0$. Replacing this in \eqref{codimLLT-2} and adding up \eqref{W2} we get \eqref{codimorbit}.
\end{proof}

The advantage of \eqref{codimorbit} against the formula provided in \cite{DeEd95} is that \eqref{codimorbit} explicitly shows the contribution of each term from the Weyr characteristics of the left and right minimal indices and the partial multiplicities of $L$ in the codimension of the orbit.

The formula \eqref{codimorbit} can be found in \cite[Prop. 3.2.2]{weyr_book} for the orbit of a matrix under similarity, which is a particular case of the orbit of a pencil under strict equivalence. In this case, only the term \eqref{W2} appears. Actually, a similar proof of the identity between \eqref{jordansum} and \eqref{W2} to the one provided here can be found in \cite[pp. 106--107]{weyr_book}.

\section{Main results}\label{main_sec}

This section contains two different kinds of results. More precisely, in Section \ref{pencils_sec} we deal with the dimension of bundles of matrix pencils, whereas Section \ref{closures_sec} deals with the closure of bundles of matrix polynomials of grade larger than $1$.

\subsection{Dimension inequalities for bundles of matrix pencils}\label{pencils_sec}

The following result, stated in \cite{Hoyo90} (based on \cite{Pokr86}), is key to prove the main result of this section (Theorem \ref{T16}).

\begin{theorem}\label{dehoyos_th}
If $L$ and $M$ are two matrix pencils of the same size and $h:=\rank L-\rank M\geq0$, then $M\in\overline\orb(L)$ if and only if the following three majorizations hold:
\begin{enumerate}[label=\rm{(\alph*)}]
    \item \label{war1} $(r_0(M),r_1(M),\ldots) \prec (r_0(L),r_1(L),\ldots)+(h,h,\ldots)$;
    \item \label{war2} $(\ell_0(M),\ell_1(M),\ldots) \prec (\ell_0(L),\ell_1(L),\ldots)+(h,h,\ldots)$;
    \item \label{war3} $(W_1(\lambda,L),W_2(\lambda,L),\ldots) \prec (W_1(\lambda,M),W_2(\lambda,M),\ldots)+(h,h,\ldots),$ for all $\lambda\in \overline{\CC}$.
\end{enumerate}
\end{theorem}

We also use Lemma \ref{dimineq_lem} to prove Theorem \ref{T16}. We emphasize that Lemma \ref{dimineq_lem} is the key technical result of this section and that its proof requires considerable effort.

\begin{lemma}\label{dimineq_lem}
If $L$ and $M$ are two $m\times n$ pencils such that $M\in \overline{\mathcal{O}}(L)$, then $\text{dim } \mathcal{B}(M)\leq\text{dim } \mathcal{B}(L)$. Moreover, the inequality is strict if and only if $ M\not\in \mathcal{O}(L)$.
\end{lemma}
\begin{proof}
Let $M\in \overline\orb(L)$ and set $h:=\text{rank }L-\text{rank }M$. It must be $h\geqslant 0$, by the lower semicontinuity of the rank (see, for instance, \cite[p. 283]{Hoyo90}). 

In this situation, the majorization relations in Theorem \ref{dehoyos_th} hold. Additionally, as $L$ and $M$ are of size $m\times n$,  
\begin{equation}\label{Neq}
\begin{array}{ccl}
n&=&\displaystyle\sum_{\lambda_k\in \Lambda(L)}\sum_{i=1}^\infty{W_i(\lambda_k,L)}+\sum_{i=0}^\infty{r_i(L)}+\sum_{i=0}^\infty{\ell_i}(L)-\ell_0(L)\\&=&\displaystyle\sum_{\lambda_k\in \Lambda(M)}\sum_{i=1}^\infty{W_i(\lambda_k,M)}+\sum_{i=0}^\infty{r_i(M)}+\sum_{i=0}^\infty{\ell_i}(M)-\ell_0(M).  
\end{array}
\end{equation}

By \eqref{codimbundle_def} and \eqref{codimorbit} the codimension of the corresponding bundles are 
$$
\begin{array}{ccl}
\codim \mathcal{B}(L) &=&\displaystyle\ell_0(L)n+r_0(L)m-\sum_{i=0}^\infty{r_i(L)r_{i+1}(L)}-\sum_{i=0}^\infty{\ell_i(L)\ell_{i+1}(L)}\\&&\displaystyle+\sum_{\lambda_k \in \Lambda(L)}\sum_{i=1}^\infty{W_i(\lambda_k,L)^2}-\#\Lambda(L)
\end{array}
$$ and
$$
\begin{array}{ccl}
\codim \mathcal{B}(M)&=&\displaystyle\ell_0(M)n+r_0(M)m-\sum_{i=0}^\infty{r_i(M)r_{i+1}(M)}-\sum_{i=0}^\infty{\ell_i(M)\ell_{i+1}(M)}\\&&\displaystyle+\sum_{\lambda_k \in \Lambda(M)}\sum_{i=1}^\infty{W_i(\lambda_k,M)^2}-\#\Lambda(M).
\end{array}
$$ 
Our aim is to show that $\codim \mathcal{B}(M)\geq\codim \mathcal{B}(L)$, and that $\codim\cB(L)=\codim\cB(M)$ if and only if $M\in\orb(L)$. With respect to the equality of codimensions, note that the definition of bundle implies that if $M \in \orb(L)$, then $\cB (L) = \cB (M)$ and, so, their codimensions are equal to each other. Therefore, we need to prove that for $M \in \overline\orb (L)$ the smallest $\codim \mathcal{B}(M)$ is attained only if $M \in \orb(L)$. For this purpose, let us assume that the pencil $M$ is such that
\begin{equation}\label{mincodim}
\codim \mathcal{B}(M)=\min\{\codim\cB(R):\ R\in\overline\orb(L)\}.    
\end{equation}

We carry out the proof in three steps.

\textbf{Step 1}. Let us prove that, if $\mu \in \Lambda(M) \setminus \Lambda(L)$, then there is only one Jordan block in the KCF of $M$ corresponding to $\mu$ and it has size $1\times 1$, namely $W(\mu,M)=(1,0,0,\hdots)$.

By contradiction, let us assume that there is an eigenvalue $\mu\in\Lambda(M)\setminus\Lambda(L)$ such that $W_1(\mu,M)>1$ or $W_2(\mu,M)>0$. Then 
\begin{equation}\label{tminusweyr}
    \sum_{i=1}^\infty W_i(\mu,M)^2\geq \sum_{i=1}^\infty W_i(\mu,M)> 1.
\end{equation}

Let $\widetilde M$ be a pencil whose KCF is the same as the KCF of $M$ except for the Jordan blocks corresponding to the eigenvalue $\mu$. Instead, the KCF of $\widetilde M$ has %$s$
 $s:=\sum_{i=1}^\infty W_i(\mu,M)$ new different (simple) eigenvalues, $\la_1,\hdots,\la_s$, namely $W(\la_k,\widetilde M)=(1,0,\hdots)$, for $i=1,\hdots,s$, which do not belong to $\Lambda(L)\cup\Lambda(M)$. In particular, 
\begin{equation}\label{t-s}
\#\Lambda(M)-\#\Lambda(\widetilde M)=1-s.
\end{equation}

Since $\mu,\la_k\not\in  \Lambda(L)$ and $\rank \widetilde{M} = \rank M$, the pencil $\widetilde M$ 
satisfies \ref{war1}-\ref{war3} in Theorem \ref{dehoyos_th}, so $\widetilde M\in\overline\orb(L)$.
However, from Proposition \ref{orbitcodim_prop}
$$
\begin{array}{ccl}
\codim\mathcal{B}(\widetilde M)-\codim\mathcal{B}(M)&=&\displaystyle\sum_{i=1}^\infty\sum_{k=1}^sW_i(\la_k,\widetilde M)^2-\#\Lambda(\widetilde M)\\
&&- \displaystyle\left(\sum_{i=1}^\infty{W_i(\mu,M)^2}-\#\Lambda(M)\right)\\
&=& \displaystyle s-\sum_{i=1}^\infty{W_i(\mu,M)^2}+\#\Lambda(M)-\#\Lambda(\widetilde M)\\
&=&\displaystyle s-\sum_{i=1}^\infty{W_i(\mu,M)^2}+1-s\\
&=&\displaystyle 1-\sum_{i=1}^\infty{W_i(\mu,M)^2}<0,
\end{array}
$$
where to get the third identity we have used \eqref{t-s} and for the last inequality we have used \eqref{tminusweyr}. This contradicts \eqref{mincodim}.

\textbf{Step 2}.  Let us prove that $\Lambda(M)\subseteq \Lambda(L)$. For this, we consider two cases: $h=0$ and $h\ne 0$.

\textbf{Case $h=0$:} In this case, condition \ref{war3} in Theorem \ref{dehoyos_th} implies $W_1(\lambda,L)\leq W_1(\lambda,M)$, for any $\lambda\in \Lambda(L)$. Thus, $\Lambda(L)\subseteq \Lambda(M)$.

Assume first that $M$ is regular with size $n\times n$. Then $L$ is regular with size $n\times n$ as well. 
Moreover, condition \ref{war3} in Theorem \ref{dehoyos_th} implies that $\sum_{i=1}^\infty
W_i(\lambda,L)\leq\sum_{i=1}^\infty W_i (\lambda,M)$, for any $\lambda\in\Lambda(L)$. This, combined with $\Lambda(L)\subseteq \Lambda(M)$, implies that $$n=\sum_{\lambda\in\Lambda(L)}\sum_{i=1}^\infty W_i(\la,L) \leq \sum_{\lambda\in\Lambda(L)}\sum_{i=1}^\infty W_i(\la,M) \leq \sum_{\lambda\in\Lambda(M)}\sum_{i=1}^\infty W_i(\la,M) = n.$$
Thus, all the inequalities above are equalities, which implies $\Lambda(L) = \Lambda(M)$.

Now, let us assume that $M$ is a singular pencil. 
Again, by \ref{war3} in Theorem \ref{dehoyos_th}, the inequality 
$$
\sum_{i=1}^\infty W_i(\lambda,L) \leq \sum_{i=1}^\infty W_i(\lambda,M)
$$
is satisfied for each $\lambda\in \overline{\CC}$. Next, we proceed by contradiction and assume that there exists $\mu\in\Lambda(M)\setminus \Lambda(L)$, which, by Step 1, satisfies $W(\mu,M)=(1,0,\hdots)$ and, so, 
$$
\sum_{i=1}^\infty W_i(\mu,L)=0<1= \sum_{i=1}^\infty W_i(\mu,M).
$$
Combining the two inequalities above with $\Lambda(L)\subseteq \Lambda(M)$, we conclude that
%$$\sum_{\lambda_k\in \Lambda(L)\cup\Lambda(M)}\sum_{i=1}^\infty W_i(\lambda_k,L)< \sum_{\lambda_k\in \Lambda(L)\cup\Lambda(M)}\sum_{i=1}^\infty W_i(\lambda_k,M).$$

$$\sum_{\lambda\in \Lambda(L) }\sum_{i=1}^\infty W_i(\lambda,L)< \sum_{\lambda\in \Lambda(M)}\sum_{i=1}^\infty W_i(\lambda,M).$$ 
Thus, by \eqref{Neq} and $\ell_0(M)-\ell_0(L)=h=0$, it must be
$$
\sum_{i=0}^\infty r_i(M) <\sum_{i=0}^\infty r_i(L) \qquad \textnormal{or } \qquad \sum_{i=0}^\infty \ell_i(M) <\sum_{i=0}^\infty \ell_i(L)
$$
(or both).
Without loss of generality we can assume that the first inequality holds, namely
\begin{equation}\label{assumtion_r}\sum_{i=0}^\infty r_i(M) <\sum_{i=0}^\infty r_i(L)
\end{equation}
(otherwise, we proceed in a similar way with the other inequality).

The inequality \eqref{assumtion_r} implies $r_0(L)>0$ and, then, the identity $r_0(M)-r_0(L)=h=0$ implies that $r_0(M)>0$ as well, so the KCF of $M$ has right singular blocks.

Now, let us consider a pencil $\widehat{M}$
whose KCF is the same as the one of $M$ except for the following: instead of one of the largest right
singular blocks (with size, say, $\alpha\times (\alpha+1)$) together with the $1\times 1$ Jordan block corresponding to $\mu \,\in\Lambda(M)\setminus \Lambda(L)$ in the KCF of $M$, the KCF of $\widehat M$
has a right
singular block of size $(\alpha+1)\times (\alpha+2)$ (which is the largest right singular block in the KCF of $\widehat M$). In particular, 
\begin{eqnarray}
\#\Lambda(M)-\#\Lambda(\widehat M)=1,\label{cond1}\\
W_i(\la,M)=W_i(\la,\widehat M),\quad \mbox{for all $\la\neq\mu$ and for all $i$},\label{cond2}\\
W(\mu,\widehat M)=(0,0,\hdots),\label{cond3}\\
r_i(M)=r_i(\widehat M),\quad \mbox{for $0\leq i\leq\alpha$},\label{cond4}\\
r_{\alpha+1}(M)=0,\ r_{\alpha+1}(\widehat M)=1,\label{cond5}\\
\ell_i(M)=\ell_i(\widehat M),\quad \mbox{for all $i$.}\label{cond6}
\end{eqnarray} 
Conditions \eqref{cond1}--\eqref{cond6}, together with Proposition \ref{orbitcodim_prop} and \eqref{codimbundle_def}, imply
$$
\codim \mathcal{B}(\widehat{M})-\codim \mathcal{B}(M)=-r_\alpha(\widehat M)-\#\Lambda(\widehat{M})
-W_1(\mu, M)^2+\#\Lambda(M)=-r_\alpha(\widehat M)<0.
$$

Now, let us check that conditions \ref{war1}-\ref{war3} in Theorem \ref{dehoyos_th} hold for $\widehat{M}$ instead of $M$. Condition \ref{war2} holds because of \eqref{cond6} and the fact that $M\in\overline\orb(L)$ (note that $h=0$ both for $M$ and $\widehat{M}$, since $\rank M = \rank \widehat{M}$). Condition \ref{war3} holds because of \eqref{cond2}--\eqref{cond3}, together with the fact that $M\in\overline{\orb}(L)$ and $\mu\not\in\Lambda(L)$ (again, we use that $h=0$ for $M$ and $\widehat{M}$). 

It remains to prove that condition \ref{war1} in Theorem \ref{dehoyos_th} also holds for $\widehat M$ instead of $M$. By condition \eqref{assumtion_r}, either $r_{\alpha+1}(L)>0$ or 
$r_0(M)+\cdots+r_\alpha(M)< r_0(L)+\cdots+r_\alpha(L)$. In both cases, 
$
r_0(M)+\cdots+r_j(M)\leq r_0(L)+\cdots+r_j(L)$, for all $0\leq j\leq\alpha$ (recall that $M\in\overline\orb(L)$, so condition \ref{war1} in Theorem \ref{dehoyos_th} is satisfied), and
$r_0(M)+\cdots+r_\alpha(M)+1\leq r_0(L)+\cdots+r_\alpha(L)+r_{\alpha+1}(L)
$,
and this implies $r(\widehat M)=(r_0(M),r_1(M),\hdots,r_\alpha(M),1,0,\hdots)\prec r(L)$, as wanted.

Therefore, $\widehat{M}\in \overline\orb(L)$ and $\codim \cB(\widehat{M})<\codim \cB(M)$, which is in contradiction with \eqref{mincodim}. Thus, there does not exist any $\mu\in\Lambda(M)\setminus \Lambda(L)$ and, since $\Lambda (L) \subseteq \Lambda (M)$, we get that $\Lambda (L) = \Lambda (M)$ also holds for singular pencils $M$ when $h=0$.

\medskip

\textbf{Case $h>0$:} We proceed again by contradiction and assume that there exists some $\mu\in\Lambda(M)\setminus \Lambda(L)$. Since 
$r_0(M)=r_0(L)+h>r_0(L)$, the pencil $M$ has right singular blocks. Then, we can construct a pencil $\widehat{M}$ as in the precedent case, which satisfies conditions \eqref{cond1}--\ref{cond6} and $\codim \cB(\widehat{M})<\codim \cB(M)$ as before. To check that $\widehat{M}\in \overline\orb(L)$ we only need to show that \ref{war1} in Theorem \ref{dehoyos_th} is satisfied with $\widehat{M}$ instead of $M$, since it is straightforward to see that conditions \ref{war2} and \ref{war3} hold taking into account that $h$ has the same value for $M$ and $\widehat M$ because $\rank M = \rank \widehat M$. For this, note that, from \ref{war1} in Theorem \ref{dehoyos_th} (with $L$ and $M$),
$$
\begin{array}{ccl}
r_0(M)+\cdots+r_j(M)&\leq& (r_0(L)+h)+(r_1(L)+h)+\cdots+(r_j(L)+h),
\end{array}
$$
for all $0\leq j\leq\alpha$, and
$$
\begin{array}{ccl}
r_0(M)+\cdots+r_\alpha(M)+1 &\leq& r_0(M)+r_1(M)+\cdots+r_\alpha(M)+h\\&\leq& (r_0(L)+h)+(r_1(L)+h)+\cdots+(r_\alpha(L)+h)\\&&+(r_{\alpha+1}(L)+h),
\end{array}
$$
and this implies, together with \eqref{cond4}--\eqref{cond5} again, that $r(\widehat M)=(r_0(M),\hdots,r_\alpha(M),1,0,\hdots)\prec r(L)+(h,h,\hdots)$.

Therefore, $\widehat{M}\in \overline\orb(L)$ and this contradicts \eqref{mincodim} again. So, $\Lambda (M) \subseteq \Lambda (L)$ and the proof of Step 2 ends.

\medskip

\textbf{Step 3}. We are going to see that it must be $M\in\orb(L)$. 

For this, we proceed by contradiction and assume that $M\in \overline\orb(L) \setminus \orb(L)$. From the Proposition in p. 60 of \cite{Hump75}, $\text{dim }\orb(M)<\text{dim }\mathcal{O}(L)$. Therefore, from \eqref{codimbundle_def}, 
$$
\codim \mathcal{B}(M)=2mn-\text{dim }\orb(M)-\#\Lambda(M)> 2mn-\text{dim }\orb(L)-\#\Lambda(L)=\codim \cB(L),
$$
where we have used that $\#\Lambda(M)\leq\#\Lambda(L)$, by Step 2.

Since this is in contradiction with \eqref{mincodim}, it must be $ M\in \mathcal{O}(L)$. 

Summarizing, we have proved that the smallest $\codim \cB(M)$, for $M\in \overline\orb(L)$, is realized only for $M\in\orb(L)$. For such a pencil $M$, we have $\orb(M)=\mathcal{O}(L)$, thus $\dim\cB(M)=\dim\cB(L)$.
Therefore
$\dim \cB(M)\leq\dim \cB(L),$ for all $M\in \overline\orb(L)$ and 
$\dim \cB(M)<\dim \cB(L)$ if and only if $M\in \overline\orb(L)\setminus \orb(L)$, as claimed.
\end{proof}

Now we are in the position of proving the main result in this section, which completes Theorem 16 in \cite{DD23} with the strict inequality of the dimensions of the involved bundles. 

\begin{theorem}\label{T16} Let $L$ be an $m\times n$ pencil. Then, there is a finite number of different $m\times n$ pencils, $L_1,L_2,\dots, L_p$ (with, say, $L_1=L$), satisfying $\mathcal{B}(L_i)\cap\mathcal{B}(L_j) = \emptyset$, for $i\not=j$, and such that
\begin{equation}\label{bundleunion}
\overline\cB(L)=\bigcup_{i=1}^p\cB(L_i).
\end{equation}
Moreover, $\overline\cB (L_i)$ is strictly included in $\overline\cB (L)$ and $\dim \cB(L_i)<\dim \cB(L)$, for $i\not=1$.
\end{theorem}
\begin{proof}
The existence of the decomposition \eqref{bundleunion} is proven in \cite[Th. 16]{DD23} together with the conditions that $\overline\cB (L_i)$ is strictly included in $\overline\cB (L)$ for $i\ne 1$, and $\mathcal{B}(L_i)\ne\mathcal{B}(L_j)$ for $i\ne j$. Note that the latter inequality is equivalent to $\mathcal{B}(L_i)\cap\mathcal{B}(L_j) = \emptyset$ by the definition of bundle. It remains to prove that $\dim \mathcal{B}(L_i)<\dim \mathcal{B}(L)$, for $i>1$.
    
Let us fix an index $i>1$. From the characterization of the closure of a pencil bundle given in \cite[Th. 9]{DD23} we conclude that there exists a map $\psi : \overline{\mathbb{C}} \rightarrow \overline{\mathbb{C}}$ such that $L_i \in \overline\orb(\psi_c(L))$. There are two possible cases, that we analyze separately.

    \textbf{Case 1:} $\psi$ is one-to-one over $\Lambda (L)$. In this case, it must be $L_i\in\overline\orb(\psi_c(L))\setminus \mathcal{O}(\psi_c(L))$, since otherwise $\mathcal{B}(L_i)=\mathcal{B}(\psi_c(L))=\mathcal{B}(L)$. Then, we conclude from Lemma \ref{dimineq_lem} that $\dim \mathcal{B}(L_i)<\dim \mathcal{B}(\psi_c(L))=\dim \mathcal{B}(L)$. 

    \textbf{Case 2:} $\psi$ is not one-to-one over $\Lambda (L)$. Note first that $\dim \mathcal{O}(L)=\dim\mathcal{O}(\psi_c( L))$. Indeed, the Weyr characteristics of the left and right minimal indices of $L$ coincide, respectively, with those of $\psi_c(L)$. Moreover, the sets of all nonzero elements in the union of the Weyr characteristics of all eigenvalues (counted with their multiplicities) of $L$ and $\psi_c(L)$ are equal to each other as well. This is because the Weyr characteristics of the eigenvalues of $\psi_c(L)$ that have coalesced are the union of the Weyr characteristics of the corresponding eigenvalues in $L$ (see Definition \ref{coalescence_def}). Thus the identity $\dim \mathcal{O}(L)=\dim\mathcal{O}(\psi_c( L))$ follows from Proposition \ref{orbitcodim_prop}. Moreover, since $\psi$ is not one-to-one over $\Lambda (L)$, the number of different eigenvalues of $L$ is greater than that of $\psi_c(L)$. Therefore, $\dim \mathcal{B}(\psi_c(L))<\dim \mathcal{B}(L)$. Since, by Lemma \ref{dimineq_lem}, $\dim \cB(L_i)\leq\dim \cB (\psi_c(L))$, we conclude that $\dim\cB(L_i)<\dim\cB(L)$.
\end{proof}

\subsection{Bundle closures of matrix polynomials of grade larger than $1$}\label{closures_sec}

The first main result of this part is Theorem \ref{main_th}, which extends Theorem 12 in \cite{DD23}, valid only for matrix pencils. We include here this result for completeness and for the ease of reading.

\begin{theorem}{\rm (\cite[Th. 12]{DD23}).}\label{dDD23_th} 
 Let $L$ and $M$ be two matrix pencils of the same size. Then $\overline\cB(M)\subseteq\overline\cB(L)$ if and only if $M\in\overline\orb(\psi_c(L))$, for some map $\psi:\overline\CC\rightarrow\overline\CC$.   
\end{theorem}

The following technical lemma will be key to prove Theorem \ref{main_th}.

\begin{lemma}\label{companion-inclusion_lemma} Let $P, Q \in {\rm POL}^{m\times n}_d$. Then, $\overline{\cB^{\rm syl}}(\cC_P)\subseteq\overline{\cB^{\rm syl}}(\cC_Q)$ if and only if $\cC_P\in\overline{\orb^{\rm syl}}(\cC_{\psi_c(Q)})$, for some map $\psi:\overline\CC\rightarrow\overline\CC$.
\end{lemma} 
\begin{proof} Let us first assume that $\overline{\cB^{\rm syl}}(\cC_P)\subseteq\overline{\cB^{\rm syl}}(\cC_Q)$. Since $\cC_P \in \overline{\cB^{\rm syl}}(\cC_P)$, we get that $\cC_P \in \overline{\cB^{\rm syl}}(\cC_Q) =\overline\cB(\cC_Q)\cap{\rm GSYL}^{m\times n}_d$, by Lemma \ref{sylvclosure_lem}. Thus, $\cC_P \in \overline\cB(\cC_Q)$. Combining Lemma \ref{pinclosureq_lemma} and Theorem \ref{dDD23_th}, we get $\cC_P\in\overline{\orb}(\psi_c(\cC_{Q}))$, for some map $\psi:\overline\CC\rightarrow\overline\CC$. This implies, using \eqref{psi-comp}, that $\cC_P\in\overline{\orb}(\cC_{\psi_c(Q)})$ and, since $\cC_P\in {\rm GSYL}^{m\times n}_d$, $\cC_P\in \overline{\orb}(\cC_{\psi_c(Q)}) \cap {\rm GSYL}^{m\times n}_d = \overline{\orb^{\rm syl}}(\cC_{\psi_c(Q)})$, where the last equality follows from \cite[Lemma 3.1]{DmDo17}.

Now, let us assume that $\cC_P\in \overline{\orb^{\rm syl}}(\cC_{\psi_c(Q)})$. Again by \cite[Lemma 3.1]{DmDo17}, this implies $\cC_P\in \overline{\orb}(\cC_{\psi_c(Q)}) \cap {\rm GSYL}^{m\times n}_d$. Thus, $\cC_P\in \overline{\orb}(\cC_{\psi_c(Q)}) = \overline{\orb}(\psi_c(\cC_{Q}))$, where the last equality follows from \eqref{psi-comp}. Theorem \ref{dDD23_th} implies $\overline{\cB}(\cC_{P}) \subseteq \overline{\cB}(\cC_{Q})$. Therefore,
$\overline{\cB}(\cC_{P}) \cap  {\rm GSYL}^{m\times n}_d \subseteq \overline{\cB}(\cC_{Q}) \cap {\rm GSYL}^{m\times n}_d$ and, by Lemma \ref{sylvclosure_lem}, $\overline{\cB^{\rm syl}}(\cC_P)\subseteq\overline{\cB^{\rm syl}}(\cC_Q)$.
\end{proof}

Now we are in the position to prove the first main result of this part.

\begin{theorem}\label{main_th}
   Let $P,Q \in {\rm POL}^{m\times n}_d$. Then, $\overline\cB(P)\subseteq\overline\cB(Q)$ if and only if $P\in\overline\orb(\psi_c(Q))$, for some map $\psi:\overline\CC\rightarrow\overline\CC$.
\end{theorem}
\begin{proof} By \eqref{bunequivalence}, $\overline\cB(P)\subseteq\overline\cB(Q)$ is equivalent to $\overline{\cB^{\rm syl}}(\cC_P)\subseteq\overline{\cB^{\rm syl}}(\cC_Q)$ and, by Lemma \ref{companion-inclusion_lemma}, this is in turn equivalent to $\cC_P\in\overline{\orb^{\rm syl}}(\cC_{\psi_c (Q)})$, for some map $\psi:\overline\CC\rightarrow\overline\CC$. Applying the homeomorphism \eqref{homeo}, this is also equivalent to $P=f^{-1}(\cC_P)\in f^{-1}(\overline{\orb^{\rm syl}}(\cC_{\psi_c(Q)}))$. But, by \eqref{invforb}, $f^{-1}(\overline{\orb^{\rm syl}}(\cC_{\psi_c(Q)}))=\overline\orb(\psi_c(Q))$, and this proves the statement.
\end{proof}

In plain words, what Theorem \ref{main_th}, combined with Lemma \ref{pinclosureq_lemma}, says is that a matrix polynomial $P$ is in the closure of the bundle of another matrix polynomial $Q$ if and only if $P$ belongs to the closure of the orbit of a matrix polynomial obtained from $Q$ after coalescing some eigenvalues. Therefore, the rules for the inclusion between bundle closures of matrix polynomials in ${\rm POL}^{m\times n}_d$ are those for the inclusion between orbit closures plus coalescence of eigenvalues. 

Rules for the relation of ``covering" (namely, inclusion without any intermediate orbit or bundle closure) between closures have been provided in \cite{DJKV20} for both orbits (in \cite[Th. 8]{DJKV20}) and bundles (in \cite[Th. 10]{DJKV20}) of matrix polynomials. Applying recursively these rules, we could obtain the rules for the inclusion of closures of either orbits or bundles of matrix polynomials, but it is not obvious how to obtain Theorem \ref{main_th} from this approach. 

The second main result in this part is Theorem \ref{union_th}, which is the extension of Theorem \ref{T16} from matrix pencils to matrix polynomials of arbitrary grade.

\begin{theorem}\label{union_th}
  Let $P(\la)$ be an $m\times n$ matrix polynomial of grade $d$. Then, there is a finite number of different $m\times n$ matrix polynomials of grade $d$,  $P_1(\la), P_2(\la), \hdots,\allowbreak P_\ell(\la)$ (with, say, $P_1=P$), satisfying $\cB(P_i)\cap\cB(P_j) = \emptyset$, for $i\neq j$, and such that
  \begin{equation}\label{unionbun}
      \overline\cB(P)=\bigcup_{i=1}^\ell\cB(P_i).
  \end{equation}
Moreover, $\overline\cB(P_i)$ is strictly included in $\overline\cB(P)$ and $\dim\cB(P_i)<\dim\cB(P)$, for $i\neq1$.  
\end{theorem}
\begin{proof}
By Theorem \ref{T16}, which is the particular case of the statement for matrix pencils, $\overline\cB(\cC_P)$ is the union of a finite number of bundles of matrix pencils with the same size as $\cC_P$, including $\cC_P$ itself, namely 
\begin{equation}\label{union}
\overline{\cB}(\cC_P)=\bigcup_{i=1}^s\cB(L_i),
\end{equation}
 with, say, $L_1=\cC_P$ and, by Theorem \ref{T16},
 \begin{equation}\label{strictinclusion}
 \overline{\cB}(L_i)\varsubsetneqq \overline\cB(\cC_P),
 \end{equation}
 for $i\neq1$. This strict inclusion will be used later.
 
 Replacing \eqref{union} in \eqref{fclosure} we obtain
 \begin{equation}\label{fclosure2}
 f(\overline{\cB}(P))=\left(\bigcup_{i=1}^s\cB(L_i)\right)\cap{\rm GSYL}_d^{m\times n}=\bigcup_{i=1}^s\left(\cB(L_i)\cap{\rm GSYL}_d^{m\times n}\right).
 \end{equation}
 Now, for each $1\leq i\leq s$ such that $\cB(L_i)\cap{\rm GSYL}_d^{m\times n}$ is not empty, there is at least one companion pencil, say $\cC_{P_i}$, for some $P_i\in{\rm POL}_d^{m\times n}$, such that $\cC_{P_i}\in\cB(L_i)$, so $\cB(L_i)=\cB(\cC_{P_i})$. Note that there is at least one such index, namely $i=1$, since $L_1=\cC_P$. Therefore, we may assume that the first $\ell$ bundles $\cB(L_i)$ in \eqref{fclosure2}, for some $1\leq \ell\leq s$, are of the form $\cB(\cC_{P_i})$. Then \eqref{fclosure2} becomes
 $$
 f(\overline\cB(P))=\bigcup_{i=1}^\ell\left(\cB(\cC_{P_i})\cap{\rm GSYL}_d^{m\times n}\right)=\bigcup_{i=1}^\ell \cB^{\rm syl}(\cC_{P_i}).
 $$
 This is equivalent to
 $$
 \overline\cB(P)=f^{-1}\left(\bigcup_{i=1}^\ell \cB^{\rm syl}(\cC_{P_i})\right)=\bigcup_{i=1}^\ell f^{-1}(\cB^{\rm syl}(\cC_{P_i}))=\bigcup_{i=1}^\ell \cB(P_i),
 $$
 and this concludes the proof of \eqref{unionbun}.

 Let us note that, by Lemma \ref{pinclosureq_lemma},
 \[
     \overline{\cB}(P_i)\subseteq\overline\cB(P),\qquad\mbox{for all $i=1,\hdots,\ell$.}
 \]
To prove that $\overline\cB(P_i)$ is strictly included in $\overline\cB(P)$, for $i\neq1$, assume, by contradiction, that $\overline \cB(P)=\overline\cB(P_i)$, for some $i\neq1$. This implies that $f(\overline{\cB}(P))=f(\overline{\cB}(P_i))$, namely $\overline\cB(\cC_P)\cap {\rm GSYL}_d^{m\times n}=\overline\cB(\cC_{P_i})\cap{\rm GSYL}_d^{m\times n}$, by \eqref{fclosure}. Since $\cC_P\in\overline\cB(\cC_P)\cap{\rm GSYL}_d^{m\times n}$, this implies that $\cC_P\in\overline\cB(\cC_{P_i})\cap{\rm GSYL}_d^{m\times n}$, so, in particular, $\cC_P\in\overline\cB(\cC_{P_i})$. Now, by \cite[Lemma 11]{DD23} (or, by Lemma \ref{pinclosureq_lemma} applied to pencils), we conclude that $\overline{\cB}(\cC_P)\subseteq\overline{\cB}(\cC_{P_i})$, which is a contradiction with \eqref{strictinclusion}.

It remains to prove that $\dim\cB(P_i)<\dim\cB(P)$ for $i\ne 1$. By Theorem \ref{T16}, $\dim \cB(\cC_{P_i})<\dim \cB(\cC_P)$, for $i\not=1$, namely $\codim\cB(\cC_P)<\codim\cB(\cC_{P_i})$. Additionally, from \eqref{codimbundle_def}, \eqref{codimpoly}, and \eqref{codimbcp}, we get 
$$
 \begin{array}{ccl} 
  \codim \cB(\cC_P)&=&\codim \orb(\cC_P)-\#\Lambda(\cC_P) 
\\
  &=&\codim \orb^{\rm syl}(\cC_P)-\#\Lambda(\cC_P) \\& =& \codim \cB^{\rm syl}(\cC_P)=\codim \cB(P).
\end{array}
$$
Similarly, $\codim \cB(\cC_{P_i})=\codim \cB(P_i)$, so we conclude that $\codim\cB(P)<\codim\cB(P_i)$, namely $\dim \cB(P_i)<\dim \cB(P)$, as wanted.
\end{proof}

We conclude this section with an example that illustrates results on closures of bundles.
\begin{example}
Let us consider the set of $2\times 2$ quadratic matrix polynomials (namely $m=n=d=2$). This set is the union of $19$ bundles, which are indicated in Table {\rm\ref{bundle_table}}. These bundles can be obtained from Theorem {\rm3.3} in {\rm\cite{ddvd15}}, that provides necessary and sufficient conditions for given lists to be realizable as the Weyr characteristics of the left and right minimal indices and of the partial mutliplicities of a matrix polynomial of a given rank and grade. In the first three columns of the table we indicate the Weyr characteristics of the right, left, and partial multiplicities of each matrix polynomial belonging to the corresponding bundle, and in the fourth column we display one matrix polynomial in the bundle. Finally, in the last column we show the codimension of the bundle, according to \eqref{codimbundle_def}, \eqref{codimpoly}, and \eqref{codimorbit}. We have arranged them in non-decreasing order of the codimension.

From Table  {\rm \ref{bundle_table}} and Theorem {\rm\ref{main_th}} we can identify all inclusion relationships between the closures of bundles. For instance, according to Theorem {\rm\ref{main_th}} (recall also Lemma {\rm\ref{pinclosureq_lemma})}, $P_{9}\in\overline\cB(P_3)$, because $P_{9}\in\overline\orb(\psi_c(P_3))$, with $\psi:\CC\rightarrow\CC$ being such that $\mu:=\psi(\mu_1)=\psi(\mu_2)$, for $\mu_1\neq\mu_2$. Note that $W(\mu,\psi_c(P_3))=W(\mu_1,P_3)\cup W(\mu_2,P_3)=(1,1,1,1)$, so actually $P_9\in\orb(\psi_c(P_3))$.
Another way to see that $P_{9}\in\overline\cB(P_3)$ is by finding a sequence of pencils in $\cB(P_3)$ that converges to a pencil in $\cB(P_{9})$. For instance, we may consider the sequence
$$
  \left\{ Q_n:=\begin{bmatrix}
 (\la-\mu_1)^2&1\\
    0&(\la-\mu_1)(\la-\mu_1+1/n)
\end{bmatrix}\right\}_{n\in\mathbb N},
$$
so that $Q_n\in\cB(P_3)$, for all $n\in\mathbb N$, and $Q_n$ converges to $P_{9}$ as $n\rightarrow\infty$.

Similarly, $P_3\in\overline\cB(P_2)$, since $W(\mu_1,P_3)=(1,1,1)=W(\mu_1,P_2)\cup W(\mu_3,P_2)$ (we omit the argument through the appropriate map $\psi$). Again, we can easily find a sequence of pencils in $\cB(P_2)$ converging to $P_3$, namely
$$
\left\{ R_n:=\begin{bmatrix}
 (\la-\mu_1)^2&1\\
    0&(\la-\mu_1+1/n)(\la-\mu_3)
\end{bmatrix}\right\}_{n\in\mathbb N}.
$$
To see that $R_n\in\cB(P_2)$, for all $n\in\mathbb N$ (except, perhaps at most for one value of $n$ if $\mu_3 = \mu_1 - 1/n$ for one $n$), note that $\det R_n=(\la-\mu_1)^2(\la-\mu_1+1/n)(\la-\mu_3)$ and the gcd of all $1\times 1$ minors of $R_n$ is $1$, so $W(\mu_1,R_n)=(1,1)$ and $W(\mu_1-1/n,R_n)=W(\mu_3,R_n)=(1)$.

Note also that, according to Theorem {\rm \ref{main_th}}, also $P_8\in\overline{\cB}(P_2)$, since $(1,1)\prec(2)$. In this case, we can also find a sequence in $\cB(P_2)$ which converges to $P_8$, namely: 
$$
\left\{ S_n:=\begin{bmatrix}
 (\la-\mu_1)(\la-\mu_2)&1/n\\
    0&(\la-\mu_1)(\la-\mu_3)
\end{bmatrix}\right\}_{n\in\mathbb N}.
$$
To see that $S_n\in\cB(P_2)$, for all $n\in\mathbb N$, we can argue as in the previous case, namely $\det S_n=(\la-\mu_1)^2 (\la-\mu_2) (\la-\mu_3)$ and the gcd of all $1\times1$ minors of $S_n$ is $1$, so $W(\mu_1,S_n)=(1,1)$ and $W(\mu_2,S_n)=W(\mu_3,S_n)=(1)$.

However, $P_8\not\in\overline\cB(P_3)$, since there is no way to coalesce eigenvalues in $P_3$ and apply the majorization rules \ref{war1}--\ref{war3} in Theorem {\rm\ref{dehoyos_th}}.

\begin{table}[]\footnotesize
    \centering
    \begin{tabular}{c|c|c|c|c}
         $r$&$\ell$&$W$&$P$&Codim\\\hline\hline
         $(0)$&$(0)$&$(1);(1);(1);(1)$&\footnotesize$\begin{array}{c}P_1=\begin{bmatrix}
             (\la-\mu_1)(\la-\mu_2)&0\\0&(\la-\mu_3)(\la-\mu_4)
         \end{bmatrix}\\\mu_i\neq\mu_j,\ (i\neq j)\end{array}$&$0$\\\hline
         $(0)$&$(0)$&$(1,1);(1);(1)$&\footnotesize$\begin{array}{c}P_2=\begin{bmatrix}
             (\la-\mu_1)^2&0\\0&(\la-\mu_2)(\la-\mu_3)
         \end{bmatrix}\\\mu_i\neq\mu_j,\ (i\neq j)\end{array}$&$1$\\\hline
         $(0)$&$(0)$&$(1,1,1);(1)$&\footnotesize$\begin{array}{c}P_3=\begin{bmatrix}
             (\la-\mu_1)^2&1\\0&(\la-\mu_1)(\la-\mu_2)
         \end{bmatrix}\\\mu_1\neq\mu_2\end{array}$&$2$\\\hline
         $(0)$&$(0)$&$(1,1);(1,1)$&\footnotesize$P_4=\begin{array}{c}\begin{bmatrix}
             (\la-\mu_1)^2&0\\0&(\la-\mu_2)^2
         \end{bmatrix}\\\mu_1\neq\mu_2\end{array}$&$2$\\\hline
         $(1)$&$(1,1,1)$&$(0)$&$P_{5}=\begin{bmatrix}
             1&0\\\la^2&0
         \end{bmatrix}$&$2$\\\hline
         $(1,1)$&$(1,1)$&$(0)$&$P_{6}=\begin{bmatrix}
             1&\la\\\la&\la^2
         \end{bmatrix}$&$2$\\\hline
         $(1,1,1)$&$(1)$&$(0)$&$P_{7}=\begin{bmatrix}
             1&\la^2\\0&0
         \end{bmatrix}$&$2$\\\hline
         $(0)$&$(0)$&$(2);(1);(1)$&\footnotesize$\begin{array}{c}P_8=\begin{bmatrix}
             (\la-\mu_1)(\la-\mu_2)&0\\0&(\la-\mu_1)(\la-\mu_3)
         \end{bmatrix}\\\mu_i\neq\mu_j,\ (i\neq j)\end{array}$&$3$\\\hline
         $(0)$&$(0)$&$(1,1,1,1)$&$P_9=\begin{bmatrix}
             (\la-\mu)^2&1\\0&(\la-\mu)^2
         \end{bmatrix}$&$3$\\\hline 
          $(1)$&$(1,1)$&$(1)$&$P_{10}=\begin{bmatrix}
             (\la-\mu)&0\\(\la-\mu)^2&0
         \end{bmatrix}$&$3$\\\hline
         $(1,1)$&$(1)$&$(1)$&$P_{11}=\begin{bmatrix}
             \la-\mu&(\la-\mu)^2\\0&0
         \end{bmatrix}$&$3$
        \\\hline
         $(1)$&$(1)$&$(1);(1)$&$\begin{array}{c}P_{12}=\begin{bmatrix}
             (\la-\mu_1)(\la-\mu_2)&0\\0&0
         \end{bmatrix}\\\mu_1\neq\mu_2\end{array}$&$4$  
         \\\hline
         $(0)$&$(0)$&$(2,1);(1)$&\footnotesize$\begin{array}{c}P_{13}=\begin{bmatrix}
             (\la-\mu_1)^2&0\\0&(\la-\mu_1)(\la-\mu_2)
         \end{bmatrix}\\\mu_1\neq\mu_2\end{array}$&$4$
         \\\hline
         $(0)$&$(0)$&$(1,1);(2)$&\footnotesize$\begin{array}{c}P_{14}=\begin{bmatrix}
             (\la-\mu_1)(\la-\mu_2)&\la-\mu_2\\0&(\la-\mu_1)(\la-\mu_2)\end{bmatrix}\\\mu_1\neq\mu_2\end{array}$&$4$\\\hline$(0)$&$(0)$&$(2,1,1)$&$P_{15}=\begin{bmatrix}
             (\la-\mu)^2&\la-\mu\\0&(\la-\mu)^2
         \end{bmatrix}$&$5$\\\hline
         $(1)$&$(1)$&$(1,1)$&$P_{16}=\begin{bmatrix}
             (\la-\mu)^2&0\\0&0
         \end{bmatrix}$&$5$\\\hline
         $(0)$&$(0)$&$(2);(2)$&\footnotesize$\begin{array}{c}P_{17}=\begin{bmatrix}
             (\la-\mu_1)(\la-\mu_2)&0\\0&(\la-\mu_1)(\la-\mu_2)
         \end{bmatrix}\\\mu_1\neq\mu_2\end{array}$&$6$\\\hline
         $(0)$&$(0)$&$(2,2)$&$P_{18}=\begin{bmatrix}
             (\la-\mu)^2&0\\0&(\la-\mu)^2
         \end{bmatrix}$&$7$\\\hline
         $(2)$&$(2)$&$(0)$&$P_{19}=\begin{bmatrix}
             0&0\\0&0
         \end{bmatrix}$&$8$\\\hline
    \end{tabular}
    \caption{All bundles of $2\times2$ quadratic matrix polynomials.}
    \label{bundle_table}
\end{table}
\end{example}

\section{Conclusions and open questions}\label{conclusions_sec}\color{black}
It is known that the closure of a bundle of a matrix pencil $L$ is the finite union of the bundle of $L$ together with other bundles. We have proved that the dimension of each of these other bundles is strictly smaller than the dimension of the bundle of $L$. We have also proved that the closure of the bundle of any matrix polynomial $P$ of grade larger than $1$ is the finite union of the bundle of $P$ together with other bundles of strictly smaller dimension. We have also provided a new formula for the codimension of the orbit of a matrix pencil $L$ in terms of the Weyr characteristics of the right and left minimal indices and of the partial multiplicities of $L$, and we have provided a characterization for a given matrix polynomial of grade larger than $1$ to belong to the closure of the bundle of another one, a result that extends the one for matrix pencils provided in \cite{DD23}. The extension of these results to the case of structured matrix pencils and matrix polynomials \cite{MMMM06b} remains completely open and is a source of (likely hard) open problems in this context. We also mention that codimension counts of several related problems have recently attracted the attention of different researchers (see, for instance, \cite{DeDo11,DmJK13,DmKS13,DmKS14,KaKK11}). Some of these codimension counts end up in quite complicated expressions written in terms of the Segre characteristics of the involved invariants. Another natural continuation of our work is to investigate the extent to which the use of the corresponding Weyr characteristics may lead to relevant simplifications as in the case of Proposition \ref{orbitcodim_prop}.

\smallskip

\noindent{\bf Acknowledgements}. The work of F. De Ter\'an and F. M. Dopico has been partially supported by the Agencia Estatal de Investigaci\'on of Spain through grants PID2019-106362GB-I00 MCIN/AEI/10.13039/ 501100011033/ and RED2022-134176-T, and by the Madrid Government (Comunidad de Madrid-Spain) under the Multiannual Agreement with UC3M in the line of Excellence of University Professors (EPUC3M23), and in the context of the V PRICIT (Regional Programme of Research and Technological Innovation).

The research cooperation was funded by the program Excellence Initiative - Research University at the Jagiellonian University in Krak\'{o}w.

\bibliographystyle{siamplain}
  \bibliography{library}

\end{document}